\documentclass[a4paper, a4wide, 10pt, reqno]{amsart}
\usepackage{a4wide}

\usepackage{amsmath,amssymb,amsfonts,dsfont, comment, mathtools}
\usepackage[foot]{amsaddr}
\usepackage{enumitem}
\usepackage{comment}
\usepackage{graphicx}
\usepackage{subcaption}  
\usepackage[margin=1in]{geometry} 
\usepackage{xcolor}
\theoremstyle{plain}
\newtheorem{theorem}{Theorem}
\newtheorem{definition}{Definition}
\newtheorem{proposition}{Proposition}

\usepackage{setspace}
\usepackage[authoryear]{natbib}
\setcitestyle{authoryear,open={(},close={)},comma, sort}
\renewcommand{\qed}{\hfill$\square$}
\newcommand{\rd}{\mathrm{d}}
\usepackage[hidelinks]{hyperref}

\usepackage{etoolbox}
\makeatletter

\pretocmd{\NAT@citex}{%
  \let\NAT@hyper@\NAT@hyper@citex
  \def\NAT@postnote{#2}%
  \setcounter{NAT@total@cites}{0}%
  \setcounter{NAT@count@cites}{0}%
  \forcsvlist{\stepcounter{NAT@total@cites}\@gobble}{#3}}{}{}
\newcounter{NAT@total@cites}
\newcounter{NAT@count@cites}
\def\NAT@postnote{}

\def\NAT@hyper@citex#1{%
  \stepcounter{NAT@count@cites}%
  \hyper@natlinkstart{\@citeb\@extra@b@citeb}#1%
  \ifnumequal{\value{NAT@count@cites}}{\value{NAT@total@cites}}
    {\ifNAT@swa\else\if*\NAT@postnote*\else%
     \NAT@cmt\NAT@postnote\global\def\NAT@postnote{}\fi\fi}{}%
  \ifNAT@swa\else\if\relax\NAT@date\relax
  \else\NAT@@close\global\let\NAT@nm\@empty\fi\fi
  \hyper@natlinkend}
\renewcommand\hyper@natlinkbreak[2]{#1}

\patchcmd{\NAT@citex}
  {\ifNAT@swa\else\if*#2*\else\NAT@cmt#2\fi
   \if\relax\NAT@date\relax\else\NAT@@close\fi\fi}{}{}{}
\patchcmd{\NAT@citex}
  {\if\relax\NAT@date\relax\NAT@def@citea\else\NAT@def@citea@close\fi}
  {\if\relax\NAT@date\relax\NAT@def@citea\else\NAT@def@citea@space\fi}{}{}

\makeatother

\usepackage{rotating}
\usepackage{fancyvrb}
\title[Scaling laws for linear chance-constrained programs]{Optimization under rare events: 
scaling laws for \\linear chance-constrained programs}

\author[J. Blanchet]{Jose Blanchet$^1$}
\address{$^1$Stanford University, Department of Management Science and Engineering}
\author[J. Jorritsma]{Joost Jorritsma$^2$}
\address{$^2$University of Oxford, United Kingdom}
\author[B. Zwart]{Bert Zwart$^{2,3}$}
\address{$^3$Eindhoven University of Technology, Department of Mathematics and Computer Science}
\email{jose.blanchet@stanford.edu, joost.jorritsma@stats.ox.ac.uk, bert.zwart@cwi.nl}

\begin{document}
\begin{abstract}%
We consider a class of chance-constrained programs in which profit needs to be maximized while enforcing that a given adverse event remains rare. Using techniques from large deviations and extreme value theory, we show how the optimal value scales as the prescribed bound on the violation probability becomes small and how convex programs emerge in the limit. We use our results to analyze the performance of existing popular approaches in the rare-event regime. We show that the popular CVaR and sample approximations have optimality properties under light-tailed assumptions on the randomness, while they behave sub-optimal in a heavy-tailed setting. Our results are derived using large deviations theory, extreme value theory, process techniques, and random set theory.
\end{abstract}
\maketitle
\vspace*{-0.4cm}
{\footnotesize
\hspace{1em}Keywords: chance-constrained optimization, rare events, extreme values, regular variation.

\hspace{1em}MSC Class: 90C15 (primary), 60F10, 60G70 (secondary).
}

\section{Introduction}

We analyze a class of optimization problems
where the goal is to maximize profit, subject to a reliability constraint that requires the probability of an undesirable event below a small level $\delta$. These formulations, known as chance-constrained (CC) optimization problems, arise in a wide range of applications, including:
\begin{itemize}
    \item[--] {\em Financial planning.} Banks and insurance companies need to optimize their resources such that the probability of insolvency is kept small. 
     For example, in fixed income portfolio optimization, an investment grade portfolio has a default rate of about $\delta=10^{-4}$, see \cite{frank2008municipal}. 
    \item[--] {\em Network design.} In ultra-reliable communication system design for 5G wireless networks, the failure probability should be no bigger than $\delta=10^{-5}$, as is explained in \cite{alsenwi2019chance}. 
     \item[--] {\em Power grid operations.} The European Union has set a goal that the amount of reserves available to mitigate possible frequency disturbance should be sufficient, up to once every 20 years. As these reserves are typically procured on an hourly basis, this amounts to a safety constraint of $\delta$ smaller than $10^{-5}$,
     see \cite{entso-e}.  
\end{itemize}
Additional examples originate in the {\em design of bridges} and {\em control of aircraft}, cf.\ \cite{flightcontrol}, where the probability of a rare event is sometimes required to be below $\delta = 10^{-9}$. For more examples, see \cite{Subramanyam2020, tong2020optimization}. 

Because of their direct interpretation and versatility to be applied to various important areas, 
chance-constrained optimization has been a subject of interest in Operations Research for several decades. It has become one of the standard formulations in risk-constrained optimization and is covered as such in several textbooks cf.\ \cite{prekopa1970probabilistic, prekopa2003probabilistic,  Romeijnders20, Shapiro2021}. 
Unfortunately, chance-constrained problems are typically not very tractable: the set of feasible solutions may be non-convex, and evaluating the chance constraint can be a computational challenge. In general, CC problems are provably NP-hard, see \cite{luedtke2010integer}. As a result, researchers often resort to various approaches that provide bounds or relaxations to approximate the solution to CC optimization problems. 

Our goal in this paper is to leverage the fact that in many applications, such as those described earlier, $\delta$ is very close to zero. We use this feature, together with large-deviations principles and extreme-value theory to gain insights into the asymptotic behavior (as $\delta \downarrow 0$) of a stylized family of CC problems and their relaxations. In order to provide a more detailed description of our results and their implications, let us provide a brief overview of various methods in CC optimization.

One method, proposed in \cite{nemirovski2006convex}, constructs a convex restriction of the feasible region by upper bounding the reliability constraint by the so-called conditional value at risk (abbreviated with CVaR; henceforth we refer to this as the CVaR approach, see also \cite{RockUry02}). Another option, popular because of its ease of implementation, is the so-called scenario approach, introduced by \cite{calafiore2006scenario}. The scenario approach consists of replacing the probability constraint by several Monte Carlo-based scenarios, \emph{all of which} need to satisfy the reliability constraint. If the undesirable event needs to be kept rare, the number of required samples is large, which can be prohibitive, though 
it is sometimes possible to reduce the number of samples using importance sampling techniques, see \cite{nemirovski2006scenario, barrera2016chance, BZZ23}.

Recently, several authors (see e.g.\ \cite{Subramanyam2020, NNZ19, tong2020optimization}) have suggested a different method, namely to approximate the probability constraint using ideas from large-deviations theory. The approach can be described as follows. Suppose that ${\bf P}(G(x)) \leq \delta$ is a safety constraint, $x$ is the optimization decision variable, and $G(x)$ is an undesirable event. A (heuristic) large-deviations approximation of ${\bf P} (G(x))$ can be obtained by replacing the integral of the density of an underlying probabilistic model (say $f$) over the set $G(x)$ by a supremum, i.e.: 
\begin{equation}
\label{ld-approx}
    {\bf P}(G(x)) = \int_{G(x)} f(u) \rd u \approx \sup_{u\in G(x)} f(u) = \exp\Big(- \inf_{u\in G(x)} \big\{\log 1/f(u)\big\}\Big).
\end{equation}
The previous approximation can often be made rigorous by introducing a suitable asymptotic regime. Using \eqref{ld-approx}, the chance-constrained program can be approximated by a bi-level (or nested) optimization problem (involving $x$ at the top level and $u$ at the bottom level). 

The papers \cite{Subramanyam2020}, \cite{NNZ19}, and \cite{tong2020optimization} focus on various computational aspects of this bi-level procedure. The analysis in those papers requires sufficient tractability for the underlying density $f$, and this is carried out in situations where (for example) the randomness is a mixture of Gaussian distributions or where the randomness arises from another specific light-tailed distribution. Further, this analysis focuses on developing approximations as in \eqref{ld-approx} that result in computationally tractable procedures. 

Moving on to a more precise description of our contributions, this paper takes a different perspective:
we {\em derive} rigorous 
approximations like \eqref{ld-approx} 
allowing both light and heavy-tailed distributions. In the latter case, the analogues of \eqref{ld-approx} can be harder to formulate. In turn, we consider a sequence of CC problems indexed by the safety parameter $\delta$, and complete the following technical tasks: 
\begin{enumerate}[leftmargin=17pt]
    \item[(a)] we show how large-deviations approximations for the solution of CC problems {\em emerge} in the limit and analyze the scaling of the chance-constrained program's optimal value $v_\delta$  as $\delta\downarrow 0$; 
    \item[(b)]we compare $v_\delta$ both to the performance of the CVaR relaxation and the scenario approach;  
    \item[(c)]we identify conditions for which the latter approaches are asymptotically feasible and optimal or not. 
\end{enumerate}
An attractive feature of our analysis is that it is rooted in multivariate extreme-value theory, 
allowing us to consider classes of distributions that are either light-tailed or heavy-tailed with a rich dependence structure. We restrict our analysis to a class of linear programs to keep it crisp. Specifically, the insights that we obtained based on our technical development can be summarized as follows: 
\begin{enumerate}[leftmargin=17pt]
    \item If the underlying random variables are light-tailed, then we show (Theorem~\ref{thm:lt-ccp} below) that the basic formulation of the chance-constrained program, given in \eqref{chance-constraint-opt} below, 
    converges to a convex optimization problem as $\delta\downarrow 0$, and that an appropriate version of \eqref{ld-approx} emerges. While this leads to a characterization of the growth rate of $v_\delta$, it is nontrivial to construct a sequence of asymptotically optimal solutions from the limiting optimization problem: essentially, the limit of the logarithmic estimate \eqref{ld-approx} is rather rough to guarantee feasibility for fixed $\delta>0$. However, the asymptotics are sharp enough to derive solutions with any prescribed (non-negative) asymptotic optimality gap. 
    \item Also in the context of light-tailed distributions, we show that both the CC problem and its CVaR relaxation converge to the same convex program in the regime $\delta\downarrow 0$. Thus, CVaR provides a sequence of convex approximations that are both conservative and asymptotically optimal: their optimal values scale in the same way as the actual profit. For a precise formulation, see Theorem~\ref{thm:lt-cvar}. In short, for light-tailed settings, the CVaR relaxation is asymptotically exact.
    \item  A similar, but more nuanced conclusion, once again in the light-tailed environment, holds for the scenario approach: in this case, we show that the sequence of random linear programs, indexed by the number of scenarios, converges after appropriate normalization to a deterministic optimization problem as the number of scenarios becomes large. To match the behavior of \eqref{chance-constraint-opt}, the scenario approach requires a careful choice in the number of scenarios $k$ as a function of $\delta$. The standard prescription provided by \cite{calafiore2006scenario}, see \eqref{cc-choice} below, allows for such a choice, providing a feasible and asymptotically optimal solution with probability tending to one.  

\item If the underlying random variables are heavy-tailed (more precisely, if their tails are of regular variation), we show (Theorem~\ref{thm:ht-ccp}) that our basic chance-constrained problem converges to a convex optimization problem after appropriate renormalization. The approximation \eqref{ld-approx} is not adequate here: essentially logarithmic asymptotics are not enough and a more precise asymptotic estimate is needed. We provide such an estimate  using multivariate regular variation.

\item A normalized sequence of CVaR constraints converges to a similar limit as the chance constraints. However, CVaR yields a value of the objective function which is a constant factor away from the actual optimum as $\delta\downarrow0$. 
This makes CVaR asymptotically suboptimal in the heavy-tailed case.
A detailed statement is Theorem~\ref{thm:ht-cvar} below. 
    \item The scenario approach behaves poorly in the heavy-tailed case: the appropriately normalized sequence of random optimization problems converges to a limit, but this limit is still random and feasibility to the original problem cannot be guaranteed, cf.\ Theorem~\ref{thm:ht-scen} below. This is in contrast to methods such as sample-average approximation (also known) as empirical optimization in which, as the Monte Carlo sample size increases, the sample-based approximation converges to the desired solution. Our results suggest that in the heavy-tailed case the scenario approach must be modified to obtain asymptotically both optimal and feasible solutions. We utilize concepts from point-process theory by \cite{Resnick07} and random set theory, see \cite{Molchanov}, to establish these results. 
\end{enumerate}

The rest of this paper is organized as follows. 
We define our model and state our assumptions in Section~\ref{sec-problem}. In Section \ref{sec-light}, we derive our results for light-tailed distributions. 
Our results for heavy-tailed distributions can be found in Section \ref{sec-heavy}. In Section \ref{sec-example} we consider an example showing how it is possible to model dependence between random variables in our framework for light-tailed distributions, and how using approximations like \eqref{ld-approx} impacts feasibility. This example allows for explicit calculations. Afterwards, we verify our theoretical results with a numerical example on portfolio optimization in Section \ref{sec-numeric}. 
We give a short reflection in Section \ref{sec-conc}, and present our proofs in Section \ref{sec-proofs}, which also develops a general-purpose result that circumvents the use of epi-convergence.

\section{Problem description}
\label{sec-problem}

\noindent
{\bf Notational conventions.}
Given vectors $x,y$ of identical dimension, the inequalities $x\geq y$ and $x>y$ are defined component wise, i.e. if $x>y$ if $x_i>y_i$ for all elements $i$. For an $m$-dimensional vector $x$ we denote $|x|=\sum_{i=1}^m |x_i|$ throughout the paper.

\smallskip
\noindent{\bf Basic formulation.} Let $X$ be a compact convex subset of $[0,\infty)^m$ for some $m\in\mathbb{N}$, and assume $X$ contains the set $[0,h]^m$ for some $h>0$. Let $n\in\mathbb{N}$, and $\phi: X\times[0,\infty)^n\to \mathbb{R}$ be a function. Consider the problem 
\begin{align}
\label{chance-constraint-opt}%
\mbox{maximize }  c^Tx \hspace{1cm}
\mbox{ subject to } x\in X, \mbox{ and } {\bf P}(\phi(x,L) > 1)\le \delta,
\tag{$\rm{CCP}_\delta$}%
\end{align}
where $L \geq 0$
is an $n$-dimensional random vector. The elements
of $L$ are often referred to as risk factors and may be dependent. 
We assume $\phi$ can be written in terms of $m\times n$ matrices $A^{(i)}, i=1,...,d$ with non-negative entries: we take %
\begin{equation}
\label{phiform}
    \phi(x,L) := \max_{i=1,...,d} x^T A^{(i)} L.
\end{equation}
This formulation is general enough to model linear programs with random non-negative coefficients in the constraints. 
To see this, let $C$ be a random matrix with column vectors $C_j$, and $b$ be a random vector. We can write the matrix inequality $x^T C\leq b$ as 
$\sup_i x^T C_i/b_i  \leq 1$.
Next, define $L$ as the row vector consisting of all the $C_i/b_i$ and define the matrix $A^{(i)}$ consisting of only $0$'s and $1$'s (in particular, the entries are nonnegative) in such a way that $ A^{(i)}L = C_i/b_i$. We see also that the elements of $L$ are non-negative, if the elements of $C_i/b_i$ are nonnegative.

Our choice of $\phi$ implies that $0$ is always a feasible solution, regardless of the value of $\delta$.
In the context of linear random stochastic inequalities of the form $x^T AL \leq b$, if there exists an $x_0$ such that ${\bf P} (x_0^T AL \leq b) =1$, one can always make a change of variables by setting $z=x_0+x$ and considering $z^T AL \leq b- x_0^T AL$. In this case, $z=0$ is always feasible and our results would examine the speed of convergence of the objective function towards $c^T z_0$ as $\delta\downarrow 0$. 
In this note, we focus on the specific form (\ref{phiform}) of $\phi$ to allow for short proofs. For a more general setting, we refer to the more recent paper by \cite{deo2024characterizing}.

\smallskip
As mentioned earlier, the parameter $\delta> 0$ denotes the risk tolerance. 
We focus on the case where the reliability constraint is a {\em rare event}, i.e.\ $\delta$ is close to $0$.
We aim to understand how the optimal profit 
scales as $\delta\downarrow 0$. 

\begin{definition}[Asymptotic optimality]
Let $v_\delta$ be the optimal value of $c^Tx$ in \eqref{chance-constraint-opt}. Let $(x_\delta)$ be a sequence of $m$-dimensional vectors such that $x_\delta\in X$ that is also asymptotically feasible, i.e., 
\begin{equation}
\label{asymptotic-feasibility}
    \limsup_{\delta\downarrow 0}
     {\bf P}\big(\phi(x_\delta,L) > 1
     \big)/\delta \leq 1. 
\end{equation}
We call $(x_\delta)_{\delta>0}$ asymptotically optimal if it is asymptotically feasible and
\begin{equation}
    \label{asymptotic-optimality}
    c^T x_\delta\, \big/\, v_\delta\rightarrow 1,\qquad \text{as }\delta\downarrow 0.
\end{equation}  
\end{definition}

We now introduce two formulations that are  
popular alternatives to \eqref{chance-constraint-opt}; we are interested in whether their solutions are feasible and asymptotically optimal. The determination of asymptotic optimality is facilitated by
Theorems~\ref{thm:lt-ccp} and~\ref{thm:ht-ccp} below, as they establish the asymptotic behavior of the optimal value $v_\delta$ of \eqref{chance-constraint-opt} for light-tail and heavy-tail distributions, respectively.

\medskip 
\noindent{\bf CVaR formulation.} The problem \eqref{chance-constraint-opt} is generally not convex. 
For this reason,  
 \cite{nemirovski2006convex} 
introduce a convex relaxation which typically (in particular in the case of non-atomic distributions) coincides with the {\em conditional value-at-risk} 
\begin{equation}
\label{cvar}
{\bf CVaR} (\phi(x,L))=    \inf_{\tau \in {\mathbb{R}}} \Big\{\tau + \frac 1\delta {\bf E}\big[(\phi(x,L)-\tau)^+\big]\Big\},
\end{equation}
where $x^+:=\max(0, x)$, see also \cite{RockUry02}.

In particular, \cite{nemirovski2006convex} propose to solve the optimization problem 
\begin{align}
\label{cvar-opt}%
\mbox{maximize }  c^Tx \hspace{1cm}
\mbox{ subject to } x\in X, \mbox{ and } 
{\bf CVaR} (\phi(x,L)) \leq 1.
\tag{$\rm{CVaR}_\delta$}%
\end{align}
Let $v_{{\bf CVaR}, \delta}$ be the optimal solution of the CVaR version of \eqref{chance-constraint-opt}.
The CVaR constraint in \eqref{cvar-opt} is more conservative than the chance constraint in \eqref{chance-constraint-opt}, so 
$v_{{\bf CVaR}, \delta} \leq v_\delta$.
However, as this optimization problem is convex, it may be more tractable than \eqref{chance-constraint-opt}. 
Theorems~\ref{thm:lt-cvar} and~\ref{thm:ht-cvar} below show how 
$v_{{\bf CVaR}, \delta}$
compares to $v_\delta$ as $\delta\downarrow 0$. 
\medskip

\noindent{\bf Scenario formulation.}
Another approach towards solving chance-constrained programs is to  resort to Monte-Carlo methods, where the chance constraint is
replaced by constraints of the form $\phi(x,L^{({i})}) \leq 1$, with $L^{(1)},...,L^{(k)}$ a number of samples called scenarios. 
This results in the optimization problem
\begin{align}
\label{sample-problem}%
\mbox{maximize }  c^Tx \hspace{1cm}
\mbox{subject to }  x\in X \mbox{ and } \phi(x,L^{(i)})\leq 1, \quad \forall i\in[k],
\tag{$\mathrm{SCP}_k$}%
\end{align}
where $[k]:=\{1,\ldots,k\}$.
To ensure feasibility of \eqref{chance-constraint-opt} with high enough probability, $k$ needs to be chosen appropriately as $\delta\downarrow 0$. In \cite{calafiore2006scenario} it is shown that the choice 
\begin{equation}
\label{cc-choice}
    k(\delta, \varepsilon) = \frac 2\delta \log \frac1\varepsilon + 2m+ \frac{2m}{\delta} \log \frac 2\delta
\end{equation}
gives a solution that is feasible with probability at least $1-\varepsilon$. This bound is valid for any function $\phi(x, L)$  that is convex in $x$.
We show in Theorem~\ref{thm:lt-scen}  below that the choice $k=k(\delta, \varepsilon)$ matches the behavior of the optimal value in the light-tailed case if $\varepsilon$ does not tend to 0 with $\delta$ too quickly.
\medskip 

\noindent{\bf Plan.} We analyze the above three families of optimization problems in the next two sections by making appropriate assumptions on $L$. We will apply concepts and techniques from multivariate extreme-value theory and large-deviations theory. We begin with the case in which the elements of $L$ are light-tailed.

\section{Light tails}
\label{sec-light}
In this section, we analyze \eqref{chance-constraint-opt}, \eqref{cvar-opt}, and \eqref{sample-problem} for a large class of distributions guaranteeing that all moments of $L$ are finite. In particular, we assume that $Q(x):= -\log {\bf P} (L>x)$
is {\em multivariate regularly varying}: there exists a regularly varying function $q(r)=h(r)r^\beta$ for some index $\beta>0$ and slowly-varying function $h$ ($h(ar)/h(r)\rightarrow 1$ as $r\rightarrow\infty$ for any $a>0$), and a non-negative continuous function $\lambda:[0,\infty)^m\mapsto \mathbb{R}$ such that $\inf_{x\geq 0: |x|=1} \lambda (x) >0$, and 
\begin{equation}
\label{assumption-lighttails}
    Q(rx)/q(r) \rightarrow \lambda(x),  \quad \text{as }r\rightarrow\infty.
\end{equation}
The convergence as stated above is pointwise. However, it also holds uniformly on compact sets because of regular variation, see \cite{Resnick07}.
If we divide by $|x|$ and use that $h(r|x|)/h(r) \rightarrow 1$ as $r\rightarrow \infty$ for $|x|>0$, we obtain \begin{equation}\label{eq:lambda-scaling}
     \lambda(x) = \lambda(x/|x|)\times |x|^\beta.
\end{equation}
For technical simplicity, we assume that $q(r)$ is continuous and increasing so that its inverse is well-defined. We refer to \cite{deo2024characterizing} for a more general setting.

This family of distributions includes multi-variate normal distributions ($\beta=2$), exponential and gamma distributions ($\beta=1$) and  Weibull distributions with a sub-exponential tail ($\beta<1$). The dependence structure can be rather arbitrary: the main form of dependence which is relevant is to which extent large values in the sample are correlated, which is encoded in $\lambda$. We give an example in Section~\ref{sec-example}. 

\smallskip 
We apply large-deviations theory to analyze the optimal value $v_\delta$ of \eqref{chance-constraint-opt} for small $\delta$. Set 
\begin{equation}
\label{def-I}
    I(b) := \inf_{x\in[0,\infty)^m} \big\{  \lambda (x): b^Tx \geq 1 \big\} = |b|^{-\beta}I(b/|b|), \qquad b\in[0,\infty)^m,
\end{equation}
Here, $b^T$ will represent $y^TA^{(i)}$ for some $i\in\{1,\ldots, d\}$.
We state a proposition on the role of $I$ that we prove in Section~\ref{sec-proofs}. We make the change of variable $x=y/r$, and write $f\sim g$ if $f(r)/g(r)\to 1$ as $r\rightarrow\infty$, or if $f(\delta)/g(\delta)\to 1$ as $\delta\downarrow 0$.
\begin{proposition}\label{prop:properties-I}
The set $\{b: I(b) \geq 1\}$ is convex, compact, and possesses a non-empty interior. The following asymptotics hold as $r\to\infty$, 
    \begin{align}
    \label{tail bl}
    -\log {\bf P}\big( b^T L > r\big) &\sim I(b) q(r), \\
    \label{light-tail-phi}
    -\log {\bf P}\big( \max_i y^T A^{(i)} L > r\big) &\sim  \min_{i} I\big(y^T A^{(i)}\big) q(r) =: J(y) q(r).
\end{align}
\end{proposition}
Now,
 consider the following problem, closely related to~\eqref{chance-constraint-opt},
\begin{equation}
\label{limitingoptimizationproblem-lt}
\mbox{maximize } c^T y \hspace{1cm} {\mbox{ subject to }}   y\geq 0 \mbox{ and } 
J(y) \geq 1.
\tag{$\rm{LT}$}
\end{equation}
 Formulation \eqref{limitingoptimizationproblem-lt} is motivated by making the change of variable $x=y/r$ in \eqref{chance-constraint-opt} with $r=q^{-1}(\log 1/\delta)$, replacing the probability in~\eqref{chance-constraint-opt}  by the approximation \eqref{light-tail-phi}. This ensures that the chance constraint, after taking logarithms, has the desired behavior. 
The constraint $x\in X$ is no longer relevant in (\ref{limitingoptimizationproblem-lt}) as $r\rightarrow\infty$, since we assumed $[0,h]^m \subseteq X$ for some $h>0$, and $x=y/r$.

 \smallskip
 The problem \eqref{limitingoptimizationproblem-lt} is tractable: using Sion's minimax theorem, cf.\ \cite{Sion}, strong duality holds, and since $\{J(y) \geq 1\}$ has non-empty interior, \eqref{limitingoptimizationproblem-lt} can be solved by standard Lagrangian methods.  
 In Section \ref{sec-example}, we give examples for $\lambda$ for which $I$ itself has an explicit expression.

\smallskip
Our first result is that \eqref{limitingoptimizationproblem-lt} determines the correct growth rate of $v_\delta$:
\begin{theorem}\label{thm:lt-ccp}
Let $y^*$ solve \eqref{limitingoptimizationproblem-lt}. 
Then  
$v_\delta\sim c^T y^*/q^{-1}(\log 1/\delta) $ as $\delta\downarrow 0$.
\end{theorem}

To prove this result, a verification of the pointwise convergence property \eqref{light-tail-phi} is not sufficient, as it does not imply convergence of the optima to the optimum of the pointwise limiting function. To circumvent this, we prove a general statement (Proposition \ref{prop:1} below) that we also invoke in the proofs of the subsequent theorems.

\smallskip
Based on Theorem~\ref{thm:lt-ccp}, it is tempting to conjecture that $y^*/q^{-1}(\log 1/\delta)$ is asymptotically optimal, but this may not be the case, as the asymptotic estimate \eqref{light-tail-phi} is not sharp enough to guarantee the asymptotic feasibility condition \eqref{asymptotic-feasibility}, as is exemplified in Section~\ref{sec-example}. One may instead consider $(1-\eta) y^*/q^{-1}(\log 1/\delta)$, which is asymptotically feasible for any $\eta>0$ (this follows from the proof of the next theorem in Section~\ref{sec-proof2}), and accept an optimality gap $\eta$.

\smallskip
In contrast, \eqref{cvar-opt} produces an asymptotically optimal sequence of solutions. It is more conservative than \eqref{chance-constraint-opt}, and thus feasible. Moreover, it is optimal: 

\begin{theorem}\label{thm:lt-cvar}
    $v_{{\bf CVaR}, \delta}\sim v_\delta$ 
     as $\delta\downarrow 0$.
\end{theorem}
We informally describe why Theorem~\ref{thm:lt-cvar} holds. Firstly, conditional on
$\{\phi\big(y^*/q^{-1}(\log 1/\delta),L\big)>\tau\}$, the overshoot vanishes, i.e., $\phi\big(y^*/q^{-1}(\log 1/\delta),L\big)-\tau\to 0$ almost surely. 
 Secondly, we show in Section \ref{sec-proof2} that we can approximate the optimum in \eqref{cvar-opt} from below by setting $\tau$ in \eqref{cvar} close to 1. Together, these two properties imply that the constraint in the CVaR formulation \eqref{cvar-opt} behaves asymptotically the same as the probability in \eqref{chance-constraint-opt}. This makes the CVaR formulation appealing when optimizing under rare event constraints in a light-tail setting.

\smallskip 
To assess the scenario formulation \eqref{sample-problem} in the region $k\rightarrow\infty$, 
let $\big(L^{(j)}:j=1,...,k)$ be $k$ i.i.d.\ samples of $L$.
Motivated by the insights from Theorem~\ref{thm:lt-ccp},
we consider a scaled problem. We make a similar change of variables $x=y/q^{-1} (\log k)$. We exploit that $\phi( y/r, L) =\phi(y, L)/r$ and consider
\begin{align}
\label{sample-problem-scaled-lt}%
\mbox{maximize }  c^Ty \hspace{0.5cm}
\mbox{subject to }   \phi(y,L^{(i)})\leq q^{-1}(\log k), \quad \forall i\in[k]. \tag{${\rm{SCP}\textnormal -\rm{LT}}_k$}
\end{align}
\begin{theorem}\label{thm:lt-scen}
Let $Y^{(k)}$ solve \eqref{sample-problem-scaled-lt}. 
Then, as $k\to\infty$, $c^T Y^{(k)}/q^{-1}(\log k)\rightarrow c^Ty^*$, almost surely, with $y^*$ a solution of \eqref{limitingoptimizationproblem-lt}.   
In particular, when $k=k(\delta,\varepsilon)$ defined in (\ref{cc-choice}), and  $\varepsilon=\varepsilon(\delta)$ is an arbitrary sequence such that $\delta\log(1/\varepsilon)\to0$, then
$c^T Y^{(k(\delta, \varepsilon))}/ v_\delta
\to 1$ almost surely.
\end{theorem}
The condition $\delta\log(1/\varepsilon)\to 0$ allows for some $\varepsilon=\varepsilon(\delta)\to0$. For such choices, the sequence of solutions obtained from $k=k(\delta, \varepsilon)$ is feasible with probability tending to one, and asymptotically optimal.

\smallskip 
Summarizing our findings in the light-tail case, we see that the behavior of the optimal value $v_\delta$ can be characterized using the rare-event approximation \eqref{limitingoptimizationproblem-lt}. Because this approximation is based on logarithmic asymptotics \eqref{light-tail-phi}, the rescaled optimum of $y^\ast$ may not be feasible, even for small $\delta$.

Both the sampling and the CVAR approaches provide solutions which scale like the optimal value $v_\delta$, making both approaches viable in the light-tailed case.
The standard version of the sampling approach requires a large number of scenarios $k(\delta, \varepsilon)$. In recent work, it has been shown how to adapt the sampling approach and reduce $k$ by \cite{Choi2024}.

The CVaR approach, while not as tractable as \eqref{limitingoptimizationproblem-lt}, and more conservative than \eqref{chance-constraint-opt}, is asymptotically optimal in a light-tail regime, reinforcing the fact that CVaR is one of the main approaches for solving chance constraint programs. 
In the next section, we will examine to what extent these findings hold if the underlying random variables are heavy-tailed.

\section{Heavy tails}
\label{sec-heavy}

Throughout this section, we suppose that $L$ has a distribution that is of multivariate regular variation, which covers the case in which the elements
of $L$  have a power law. 
To describe this assumption in a convenient form, we use a representation which (apart from not taking logarithms) is different from the one in Section 3, building on an equivalent characterization of regular variation developed in \cite{basrak2002characterization}). Recall that $|\cdot|$ denotes the $\ell_1$ norm. 
 Our characterization guarantees that for some $\alpha>1$, and some slowly-varying function $h$,
\[\bar F(r)= {\bf P}(|L|>r) = h(r) r^{-\alpha}.
\] 
We assume that $\bar F(t)$ is continuous and strictly decreasing (it has no atoms) to reduce technicalities. As a result, the inverse $\bar F^{-1}$ is well-defined. 
Moreover, our assumption on multivariate regular variation implies that for some random vector $\Theta$ defined on the $n$-dimensional unit sphere
\begin{equation}
\label{rvmult}
    {\bf P}\big(L/|L| \in \cdot \,\big|\, |L| > r\big) \rightarrow {\bf P}(\Theta \in \cdot), \hspace{1cm} r\rightarrow\infty. 
\end{equation}
The idea behind this representation is that, if we transform  $L$ in polar coordinates, the radius and angle become independent, conditional upon the radius getting large. 
This formulation allows dependence between extreme values of risk factors. For background on multivariate regular variation we refer to \cite{basrak2002characterization, Resnick07}. While this setting includes the important example of Pareto tails, the case of log-normal tails is not included. The case of heavy-tailed Weibull tails ($\exp(-r^\alpha), \alpha \in (0,1)$) is covered by the framework of Section 3.

\smallskip
As in the light-tailed case, we first construct an appropriately normalized version of \eqref{chance-constraint-opt} in the regime $\delta\downarrow 0$. We make the change of variables $x=y/r_\delta$, with $r_\delta:=\bar F^{-1}(\delta)= \inf \{r: {\bf P}(|L|>r) \leq \delta\}$, and write the constraint in  \eqref{chance-constraint-opt} as 
${\bf P}( \phi(y, L) > r_\delta) \leq \delta.$
Our assumptions that $\bar F(t)$ is continuous and decreasing make the chance constraint equivalent to 
\begin{equation}
{\bf P}(\phi(y,L) > r_\delta)/{\bf P}(|L| > r_\delta) \leq 1.
\end{equation}
Letting $\delta\downarrow 0$,  using \eqref{rvmult}, we show in~Section \ref{sec-proof4} that
${\bf P}(\phi(y,L) > r_\delta)/{\bf P}(|L| > r_\delta) \rightarrow {\bf E}[\phi(y,\Theta)^\alpha]$, and that
\eqref{chance-constraint-opt} 
can be approximated by  
\begin{equation}
\label{limitingoptimizationproblem}
\mbox{maximize } c^Ty \hspace{1cm} {\mbox{ subject to }}   y\geq 0 \mbox{ and } {\bf E}[\phi(y,\Theta)^\alpha]\leq 1. 
\tag{$\rm{HT}$}
\end{equation}

An appealing feature of \eqref{limitingoptimizationproblem} is that the feasible region of this problem is convex, which generally is not the case for \eqref{chance-constraint-opt}. We prove the following theorem in Section \ref{sec-proof4}:

\begin{theorem}\label{thm:ht-ccp}
    If $y^*$ solves \eqref{limitingoptimizationproblem}, then $y^*/r_\delta$ is asymptotically optimal and $v_\delta \sim c^T y^*/r_\delta$ as $\delta\downarrow 0$. 
\end{theorem}

Unlike in the light-tail case, where the optimal value of \eqref{limitingoptimizationproblem-lt} could be asymptotically \emph{infeasible}, the value $y^*$ leads to an asymptotically optimal solution as asymptotic feasibility is guaranteed. 

\smallskip
We proceed to the approximations of \eqref{cvar-opt}. In contrast to the light-tailed case, in the heavy-tailed case it is no longer asymptotically optimal: 

\begin{theorem}\label{thm:ht-cvar}
    $v_{{\bf CVaR}, \delta}\sim v_\delta( 1-1/\alpha)$ as $\delta\downarrow0$. 
\end{theorem}

The reason behind the subobtimality is that the overshoot $\phi (y,L)-r_\delta$, conditioned on  $\phi (y,L)>r_\delta$, is not of smaller order than $r_\delta$ when $\phi (y,L)$ has a heavy tail. Instead, it is of the order $r_\delta \frac{\alpha}{\alpha-1}$. 

\smallskip
We proceed to the scenario-based approach. In the previous section on the light-tailed case, we observed that the approximation \eqref{sample-problem} of \eqref{chance-constraint-opt} has some optimality properties. These properties are inherited by the convergence of the space of feasible solutions, after appropriate normalization, to a deterministic set which coincides with the feasible region of the limiting optimization problem \eqref{limitingoptimizationproblem-lt}. In the heavy-tail case, this is no longer true. In contrast, a random solution emerges, making any optimality and feasibility property elusive. 

\smallskip
To see this, consider a simple example, where the goal is to maximize a scalar $x$ subject to the constraint that ${\bf} P(xL> 1) \leq \delta$, with $L$ also scalar. 
\eqref{sample-problem} can now be formulated as maximizing $x$ subject to 
$x\max_{i=1}^k L^{(i)} \leq 1$. 
If $L$ is light-tailed, 
we can normalize $\max_{i=1}^k L^{(i)}$ such that it converges to a deterministic constant (we have done this in a more general setting to obtain Theorem~\ref{thm:lt-scen}). However, in the heavy-tail case, such a normalization is not possible. In fact, 
$\max_{i=1}^k L^{(i)}/\bar F ^{-1}(1/k)$  converges in distribution to a Fr\'echet-distributed random variable, see \cite{Resnick07}, \cite{NWZ22} for background. 

\smallskip
This scalar case provides the main insights into the sample approximation method in the heavy-tailed case. For completeness, we
 describe a general technical result 
using point-process theory commonly used in extreme-value theory, cf.\ \cite{Resnick07}, as well as random set theory as developed in (\cite{Molchanov16}).
Let  $S_{n-1}^+$ be the intersection of the unit sphere and the positive orthant. 
Let $N$ be a Poisson random measure on 
$(0,\infty] \times S_{n-1}^+$ with intensity measure $\zeta(\cdot) = \nu_\alpha (\cdot) \times {\bf P}(\Theta \in \cdot)$, where $\nu_\alpha$
has density $\alpha u^{-\alpha-1}$.
This should be interpreted as follows: if $F$ is a subset of $(0,\infty] \times S_{n-1}^+$, then $N(F)$ has a Poisson distribution with parameter 
$\zeta(F)$. 
Note that $N$ has infinitely many points, but for every $\varepsilon>0$,  the number of points in $N\big((\varepsilon, \infty] \times S_{n-1}^+\big)$ is a.s.\ finite.

Under the heavy-tailed assumptions on $L$, the point process (with $\epsilon$ being the Dirac measure)
\begin{equation}
    N_k (\cdot) = \sum_{i=1}^k \epsilon _{|L^{(i)}|/\bar F^{-1}(1/k), L^{(i)}/|L^{(i)}| } (\cdot)\label{eq:pointprocess}
\end{equation}
is converging in distribution to $N$; this is a direct application of Theorem~\ref{thm:ht-scen}.2 of \cite{Resnick07}.

\smallskip
In Section~\ref{sec-proof6}, we apply this convergence result to analyze the sample approximation method in the heavy-tailed case. 
Similar to~\eqref{sample-problem-scaled-lt} we perform a change of variables $x=y/\bar F^{-1}(1/k)$, exploit the homogeneity of $\phi$ and 
define the random variable $V_{SCP}(k)$ to be the optimal value of the optimization problem 
\begin{align}
\label{sample-problem-scaled}%
\mbox{maximize }  c^Ty \hspace{1cm}
\mbox{subject to }  y\geq 0 \mbox{ and } \phi(y,L^{(i)})\leq \bar F^{-1} (1/k), \quad \forall i\in[k].\tag{${\rm{SCP}\textnormal{-}\rm{HT}}_k$}
\end{align}
In Section \ref{sec-proof6} we apply concepts from random set theory (cf.\ \cite{Molchanov}) to show 
\begin{theorem}\label{thm:ht-scen}
    $V_{SCP}(k)$ converges in distribution to the optimal value $V^*$ of the problem 
    \begin{align}
\label{sample-problem-limit}%
\mbox{maximize }  c^Ty \hspace{1cm}
\mbox{subject to }  y\geq 0 \mbox{ and } \phi(y,z)\leq 1, \quad \forall z\in N.
\end{align}
\end{theorem}
Recall that $N$ is a random point measure. The constraints $\phi(y,z)\leq 1$ in (\ref{sample-problem-limit}) 
 are required to hold for all random points $z$ in $N$; as such, there are infinitely many constraints. However, there are only finitely many points $z$ with absolute value $\varepsilon$ or larger. Thus, it is in principle possible to verify where a specific point $y$ is feasible or not: due to the homogeneity of $\phi$ one can take this $\varepsilon$ proportional to $1/|y|$. The value $V^*$ of the optimization problem is random in general, as shown above in the simplest case where $y$ and $L$ are both scalar.  

In conclusion, comparing the three main methods to (approximately) solve \eqref{chance-constraint-opt} in the heavy-tailed case, the large-deviations approximation \eqref{limitingoptimizationproblem} yields an asymptotically optimal solution, and is characterized by a convex optimization problem. In contrast, the CVaR formulation yields a solution that stays a factor $1-1/\alpha$ away from the optimum, while the scenario approach does not converge to a deterministic solution or approximation of \eqref{chance-constraint-opt} at all. 

\section{Example for light-tailed distributions}
\label{sec-example}

While the heavy-tail case directly leads to the convex program \eqref{limitingoptimizationproblem}, the light-tail case leads to a bi-level optimization problem involving the function $I$ in~\eqref{def-I}, making it less straightforward to parse.  
The example here gives insights into how the tail behavior of $L$ and its dependence structure may impact the optimal solution of \eqref{chance-constraint-opt} as $\delta\downarrow 0$. 

Throughout this section, we assume $d=1, m=n$ and consider the problem 
\begin{align}
\label{chance-constraint-opt-example}%
\mbox{maximize }  c^Tx \hspace{1cm}
\mbox{ subject to } x\geq 0, \mbox{ and } {\bf P}\Big(\sum_{i=1}^n x_i a_iL_i \geq 1\Big)\leq \delta,
\end{align}
for some constants $a_i>0, i=1,...,n$. Suppose that $L$ satisfies
\eqref{assumption-lighttails}, with  $\lambda$ 
given by
\begin{equation}
\label{lambda-example}
    \lambda(x) := \ell(x^\beta):=\ell (x_1^\beta,...,x_n^\beta):=\bigg(\sum_{i=1}^n x_i^{\beta\theta}\bigg)^{1/\theta}, \hspace{1cm} \theta\geq 1. 
\end{equation}
Here, $\ell(x)$ corresponds to the so-called 
 Gumbel-Hougaard copula (see e.g.\ \cite{Segers2012}). 
For $\theta=1$, large coordinates in $L$ occur independently; this is no longer true if $\theta>1$. For $\theta=\infty$, we take $\ell(x) = \max_i x_i$, so that large values are perfectly correlated. 

If $\gamma=\beta\theta \leq 1$, then $I(b) = \min \{
(\sum_{i=1}^n x_i^\gamma)^{1/\theta}
: b^Tx \geq 1\}$ has as optimal solution $x_{i^*} = 1/b_{i^*}$ with $i^* = \arg\max b_i$ (with arbitrary tie-breaking rule) and $x_i=0$ for $i\neq i^*$. Consequently, 
    \begin{equation}
        I(b) = |b|_\infty^{-\beta}.
    \end{equation}
Thus, if $\gamma \leq 1$, the limiting approximation \eqref{limitingoptimizationproblem-lt} of \eqref{chance-constraint-opt} simplifies to
\begin{align}\label{chance-constraint-opt-example-lt1}%
\mbox{maximize }  c^Ty \hspace{1cm}
\mbox{ subject to } y\geq 0, \mbox{ and } \max_i a_i y_i \leq 1,
\end{align}
which has solution $y_{i}^\ast=1/a_{i}$.
    
    If $\gamma=\beta\theta > 1$, then we find using the Karush-Kuhn-Tucker (KKT) conditions that the optimal $x_i$ in the definition of $I(b)$ equals $b_i^{1/(\gamma-1)}/ \sum_j b_j^{\gamma/(\gamma-1)}$, and after some simplifications we obtain
    \begin{equation}
        I(b) =\left(\sum_{i=1}^n b_i^{\gamma/(\gamma-1)}\right)^{- (\gamma-1)/\theta}.
    \end{equation}
Thus, if $\gamma>1$, 
the limiting approximation \eqref{limitingoptimizationproblem-lt} of \eqref{chance-constraint-opt} simplifies to
\begin{align}
\label{chance-constraint-opt-example-lt2}%
\mbox{maximize }  c^Ty \hspace{1cm}
\mbox{ subject to } y\geq 0, \mbox{ and } 
\sum_{i=1}^n (a_iy_i)^{\gamma/(\gamma-1)} \leq 1.
\end{align}
Using the KKT conditions we obtain the solution
\begin{equation}
    y_i^\ast = (c_i/a_i)^{\gamma-1}/\Big(\sum_{j=1}^n (c_j/a_j)^\gamma\Big).
\end{equation}

If the coordinates of $L$ are independent in the tail ($\theta=1)$, then the nature of $I(b)$ changes at the point $\beta=1$ which is precisely the point 
where the moment-generating function of $L$ becomes finite. For  $\beta<1$, the tails of $L_i$ are heavy-tailed Weibull, for which it is known that
large values of sums are caused by a single big jump, leading to the simple form of the optimal solution. 

Interestingly, the impact of the dependence structure, summarized by $\theta$, is nontrivial, in the sense that the single-big jump intuition is not relevant in our decision-making context if the tail of $L_i$ is not very heavy, and the dependence is sufficiently strong, i.e.\ if $\beta > 1/\theta$. 

We mentioned before that asymptotic optimality in the light-tailed case is not guaranteed as the solution of \eqref{limitingoptimizationproblem-lt} after normalization may not be feasible for the problem \eqref{chance-constraint-opt}. To illustrate this in more detail in the context of our example, consider the case where $\theta=1$
and $\beta < 1$. Take also $q(u)=u^\beta$. In this case, the chance 
constraint in \eqref{chance-constraint-opt-example}, taking the scaled solution
$x_i = 1/(a_i q^{-1}(\log 1/\delta))$, leads to the probability 
${\bf P}\big(\sum_{i=1}^n L_i > q^{-1}(\log 1/\delta)\big)$. 
The distribution of $L$ is sub-exponential and we have the estimate, as $\delta\downarrow 0$, 
\begin{equation}
{\bf P}\Big(\sum_{i=1}^n L_i > q^{-1}(\log 1/\delta)\Big) \sim n {\bf P}\big(L_1> q^{-1}(\log 1/\delta)\big) = n\delta,
\end{equation}
cf.\ Chapter 3 in \cite{NWZ22}. Thus, our choice leads to a violation probability of a factor $n$ larger than it should be to satisfy \eqref{chance-constraint-opt-example}.

For $\beta=1$ we obtain a sum of i.i.d.\ exponential random variables, for which  
\begin{align}
 {\bf P}\Big(\sum_{i\in[n]}L_i> q^{-1}(\log 1/\delta)\Big)&\sim{\bf P}\big(\mathrm{Gamma}(n,1)> (\log 1/\delta)\big)\\&=\int_{u=(\log 1/\delta)}^\infty\frac{1}{\Gamma(n)}x^{n-1}{\mathrm{e}}^{-x}\rd u\\ &= \frac{1}{n!}\mathrm{e}^{-(\log 1/\delta)}(\log 1/\delta)^{n}=\delta \frac{(\log 1/\delta)^n}{n!}.
\end{align} This behavior is even worse compared to $\beta<1$ as $\delta\downarrow 0$. Thus, to ensure a solution is feasible for large $r$, one may take $(1-\gamma)y^*$ for a small $\gamma$, resulting in a solution that is close ($\gamma c^T y$) to the optimal value.

\section{Example and Numerical Experiments}
\label{sec-numeric}

To complement our theoretical findings, we now perform  numerical experiments for light- and heavy-tailed settings, demonstrating how the asymptotic predictions manifest in a practical example.

\noindent{\bf Formulation.} We consider an investment portfolio optimization problem in which an investor is interested in shorting stocks that are considered overvalued. There are several websites that can be accessed to list the top 10 overvalued stocks in various industries (e.g. see tradingview, seekingalpha, etc). We do not promote any particular set of stocks to short-sell, and mention readily available sources to stress the applicability of this problem formulation. 

We consider an investment horizon, for example, one quarter of a year. The quantity $x_i$ corresponds to the number of shares that we will short-sell of the $i$-th asset/stock. As it is commonly done in practice, we relax integer constraints. 
Let $s_i$ be the price of a single share of the $i$-th  stock today and $S_i$ be the value of the $i$-th stock's share in a quarter of a year. Prices are expressed in units of monetary risk relative to the Value-at-Risk that the manager is willing to tolerate within a given confidence level. The Value-at-Risk at level $1-\delta$ is the largest amount that can be lost with probability at most $1-\delta$. For instance, if the portfolio manager is willing to tolerate at most $\$10K$ of risk (measured by the Value-at-Risk at 95\% level) and a single share of the first asset is worth $\$1K$, then $s_1=0.1$ and $\delta=0.05$.  Since we are assuming that the stock is overvalued, $\mu_i = {\bf E}(S_i) - s_i < 0$  and we write $c_i = -\mu_i$. Define $L_i = S_i -s_i$ and note that the expected return earned by shorting the stocks according to the allocation $x=(x_1,...,x_n)^T$ at the end of a quarter of a year is precisely $c^Tx$ (which is to be maximized). The amount of money that we can lose, on the other hand, is $x^TL$. So, we are interested in maximizing the expected return subject to the Value-at-Risk at level $1-\delta$ is not larger than 1. We may add a budget constraint, say $x^Ts\leq b_0$, but we will ignore this, as this constraint is not binding as $\delta \downarrow 0$. Precisely, we are interested in solving
\begin{align}\nonumber
\mbox{maximize }  c^Tx \hspace{1cm}
\mbox{ subject to } x\geq 0, \mbox{ and } {\bf P}\bigg(\sum_{i=1}^n x_i L_i \geq 1\bigg)\leq \delta,
\end{align}
In order to test the quality of our approximations, we will use a class of models for $L$ for which this problem can be efficiently solved. We consider two types of models for $L$, one light-tailed and another one is heavy-tailed, so we can test the performance of the corresponding approximations we developed for each of these two cases. We choose $n=10$ assets. For the light-tailed case, $L$ is assumed to be Gaussian, and for the heavy-tailed case $L$ is $t$-distributed with $\alpha=5$ degrees of freedom. In both cases, $L$ has the same mean and covariance parameters, which have been selected at random; the means are uniformly distributed in [-10,1], the covariance matrix was sampled from a Wishart model, with mean identity and 13 degrees of freedom (this ensures a positive definite covariance matrix). The choice of $k$ for the scenario approach is $$
k(\delta, 0.05) = \left\lceil \frac{2n \log(2/\delta)}{\delta} + \frac{2 \log(1/.05)}{\delta} + 2n   \right\rceil.
$$ 

\medskip
\noindent\textbf{Results.}
Our theoretical results are derived for non-negative $L$, which is a technical assumption imposed to facilitate our development. While we choose our models for tractability, our choice also stresses that our insights apply more broadly\footnote{The experiments were performed in Matlab 2024, for reproducibility, the random seed was selected as rng(42).}. 

\begin{figure}
  \centering
  \begin{subfigure}[b]{0.45\textwidth}
    \centering
    \includegraphics[width=\textwidth, trim=2cm 6cm 2cm 6cm, clip]{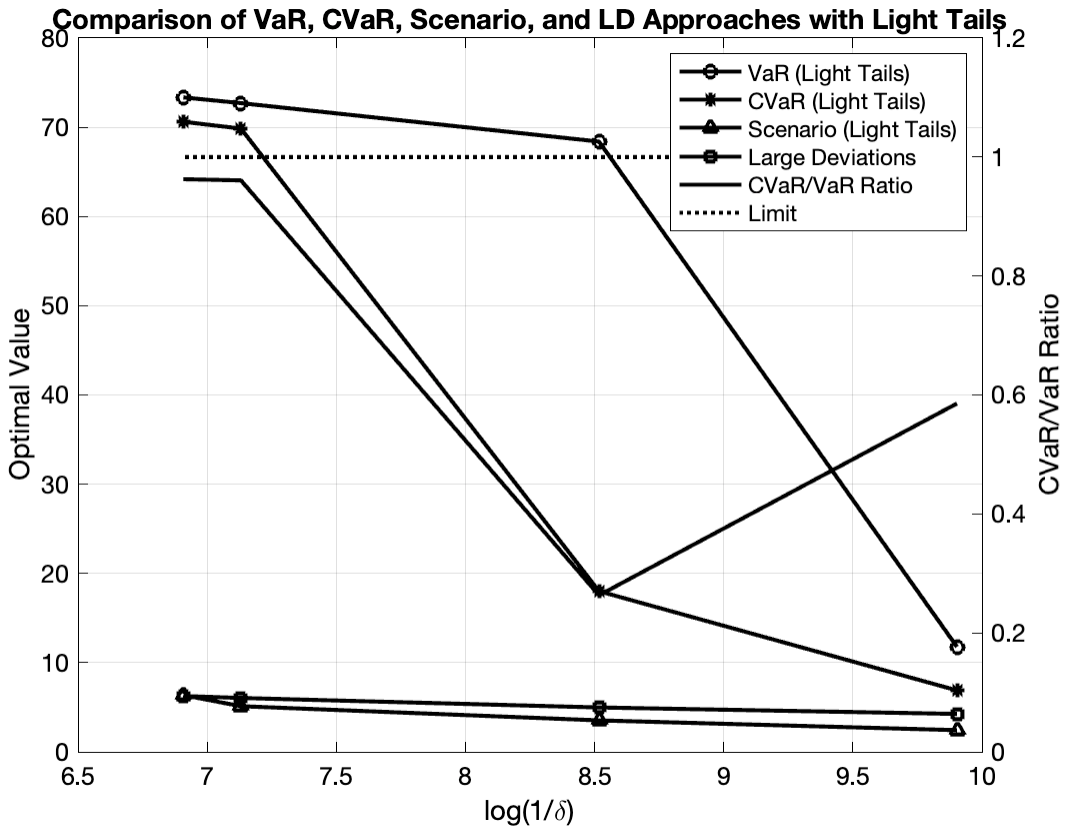}
    \caption{Light Tails - Moderate Confidence}
  \end{subfigure}
  \quad
   \begin{subfigure}[b]{0.45\textwidth}
    \centering
    \includegraphics[width=\textwidth, trim=2cm 6cm 2cm 6cm, clip]{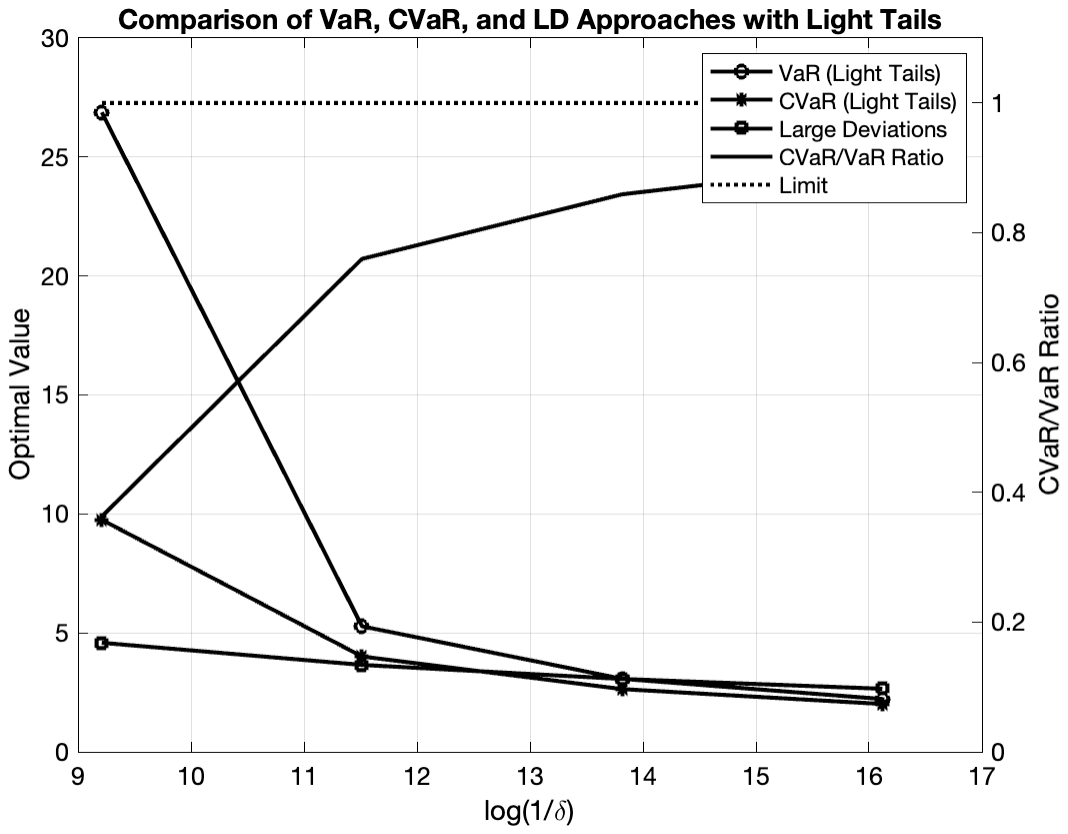}
    \caption{Light Tails - High Confidence}
  \end{subfigure}

  \vspace{1em}
   \begin{subfigure}[b]{0.45\textwidth}
    \centering
    \includegraphics[width=\textwidth, trim=2cm 6cm 2cm 6cm, clip]{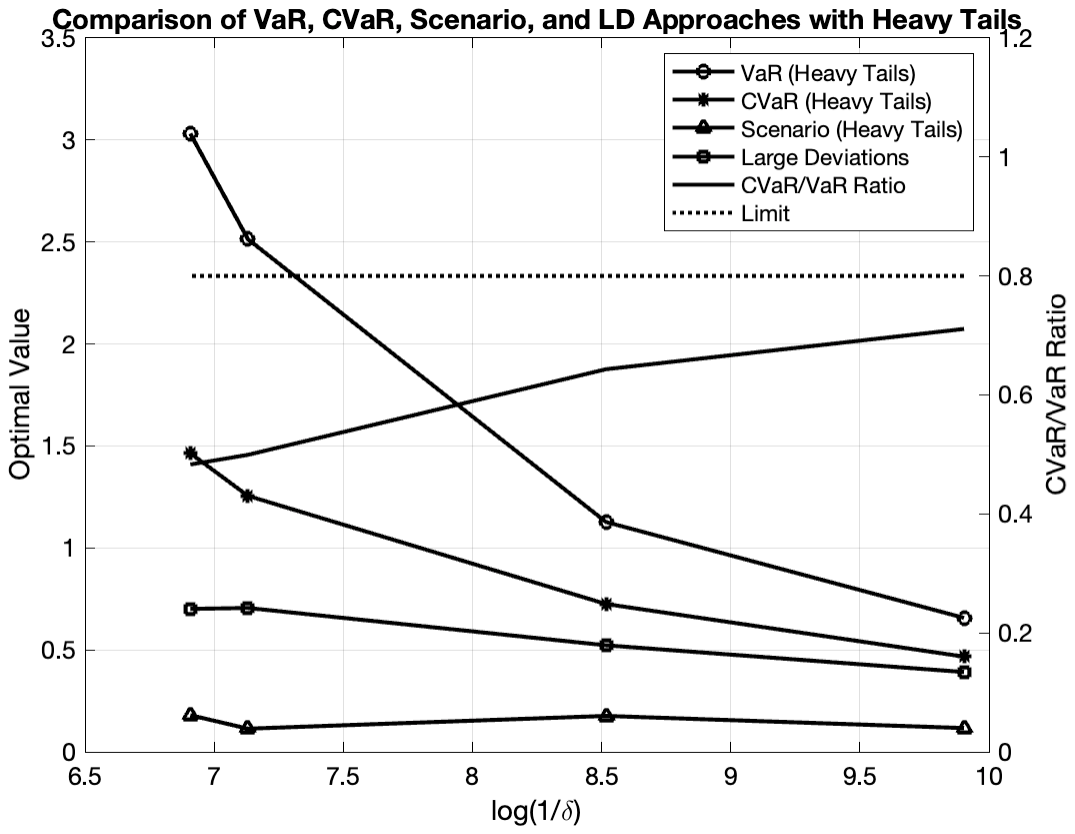}
    \caption{Heavy Tails - Moderate Confidence}
  \end{subfigure}
  \quad
  \begin{subfigure}[b]{0.45\textwidth}
    \centering
    \includegraphics[width=\textwidth, trim=2cm 6cm 2cm 6cm, clip]{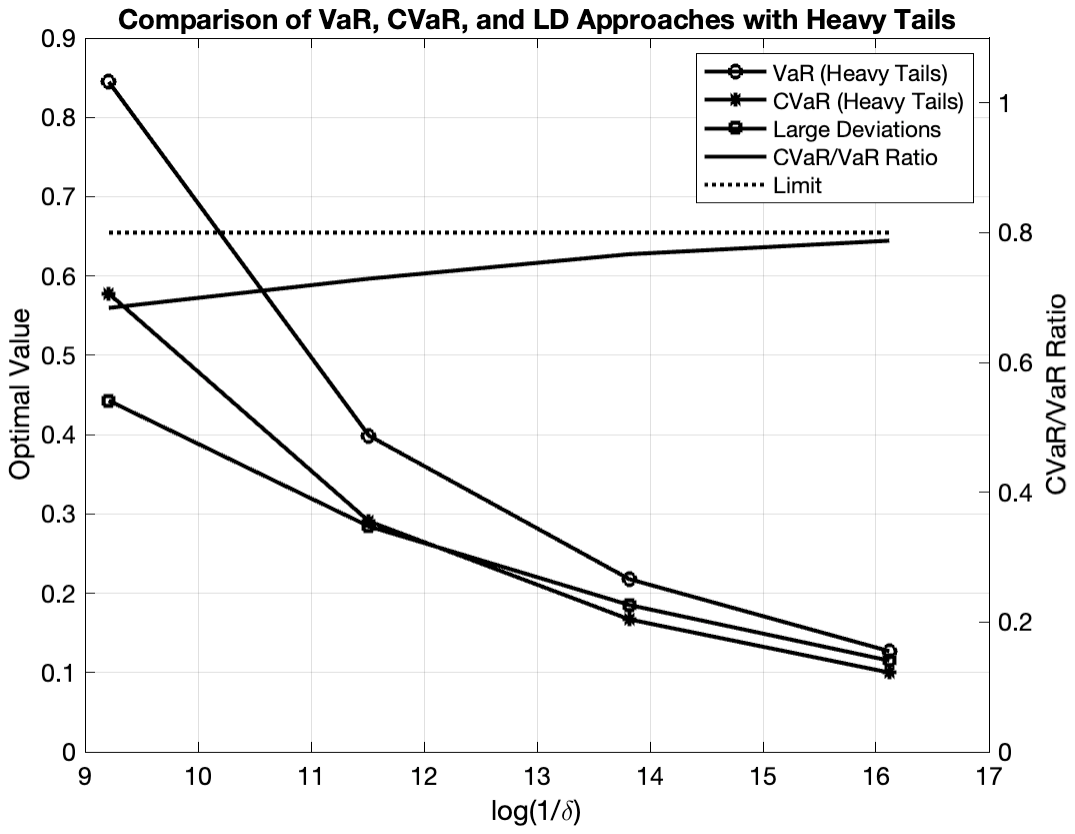}
    \caption{Heavy Tails - High Confidence}
  \end{subfigure}
    \caption{A detailed explanation of the figures is given when Figure \ref{fig:panel} is first introduced. The left vertical axis shows the exact chance constrained solution (represented as VaR) compared to the CVaR relaxation, the scenario approach and the large deviations approximation; the right vertical axis shows the ratios of CVaR/VaR solutions which converge to the theoretical limit (dotted line) as $\delta\downarrow0$. The figures on the right can be seen as a continuation of the figures on the left as the values of $\delta$ decrease.}
  \label{fig:panel}
\end{figure}

The results are shown in Figure \ref{fig:panel}, organized in two rows and two columns. The first row corresponds to the light-tailed case and the second row corresponds to the heavy-tailed case. The first column corresponds to what we call ``moderate confidence" chance constraint (enforcing values of $\delta  \in [10^{-3},5\times10^{-5}]$) and the second column corresponds to what we call ``high confidence" chance constraint (enforcing values of $\delta \in [10^{-4},10^{-7}]$). The intervals have been chosen overlapping so we can appreciate continuity in the behavior. 

Each of the four figures in the panel has two vertical axis, the left vertical axis shows the estimated optimal expected return (i.e. the optimal value of the optimization problem) under the different methods. The actual value of the chance-constraint optimization problem is denoted as VaR. We also include the CVaR approximation, the scenario approach, and the large deviations approximation. The right vertical axis plots the ratios of the CVaR to the VaR solutions. In the light-tailed case this ratio converges to 1 according to our theory, whereas in the heavy-tailed case the ratio converges to the ratio $(1-1/\alpha)=.8$ (in our case). These limiting theoretical values are represented as the ``Limit". The plots show that convergence eventually occurs as the required confidence increases (i.e. the probability of chance constrain violation decreases).

In the high confidence case, in both cases, the worst performance is obtained by the scenario approach, the large deviations approximations dominates the scenario approach.  The performance of the large deviations approximation is close to that of the scenario approach in the light tailed case, but significantly better in the heavy tailed case; the CVaR relaxation is the best performing approximation among all of the options that we consider.

In the very high confidence case, we do not present the scenario approach because it was too time consuming to produce. We can observe, as predicted by the theory, that the large deviations approximation overtakes the CVaR approximation in terms of quality of performance. The ratio VaR/CVaR is close to the predicted limit, especially for the heavy-tailed case.

The scenario approach is considered a viable and highly popular alternative, the large deviations approximations exhibits comparable performance in the light-tailed case but performs significantly better in the heavy-tailed case, with a significantly lower computational overhead. The CVaR is a viable relaxation, but eventually (in high reliability regimes) is overtaken by the large deviations approximations.

\section{Concluding comments}\label{sec-conc}
Supported by the numerical experiments, the theoretical results in this paper provide a consistent picture of how chance-constrained programs behave in the rare-event regime. We summarize our findings before proceeding to the proofs.

This paper provides the first systematic large-deviations analysis of chance-constrained optimization problems. To facilitate the development of our  algorithmic and structural insights, we focus on a linear objective and a safety constraint controlling (i.e. upper bounding by a fixed constant) the $1-\delta$ quantile of the maxima of bi-linear functions (bi-linear both in the multivariate risk factor and the decision parameter). 
We distinguish two regimes according to the tails of the risk factor, namely, light and heavy tailed regimes. 

\smallskip
Under the light-tail assumption, 
the main takeaway from Section \ref{sec-light} is that the CVaR and scenario formulations ---with carefully chosen number of scenarios--- are good approaches to chance-constrained optimization in a rare-event regime. Their solutions are asymptotically feasible (with probability tending to one), and are asymptotically optimal. The large-deviations approximation  (suitably modified as explained in Section \ref{sec-light}) can be a good approach to  chance-constrained optimization and leads to an asymptotically accurate estimate of the optimal value, though feasibility of its solution cannot be guaranteed. 

\smallskip
In the heavy-tailed regime, we explicitly show that the CVaR relaxation is no longer asymptotically optimal: its solution remains a factor $1-1/\alpha$ away from the optimum. The scenario approach behaves even worse, as explained in Section \ref{sec-heavy}. In contrast, the large-deviations approximation is asymptotically optimal in this case, as it is based on a sharper version of the large-deviations approximation  \eqref{ld-approx}. 

\smallskip
There are several natural research questions that may be considered:  (i) we expect that large-deviations approximations, for both light-tails and heavy-tails may still yield convex programming problems for functions $\phi$ which are not be bilinear; (ii)  computationally, it will be of interest to consider scaling regimes where also the dimension of the decision variable $x$ and/or the dimension of $L$ may grow ---it is kept fixed in this study; (iii) the large-deviations insights we developed could be used to develop new algorithms for chance-constrained optimization which
are asymptotically optimal in the rare event regime, and also produce close to optimal solutions for fixed values of $\delta$, just like the notion of asymptotic efficiency has led to novel algorithms in the area of rare-event simulation.

\section{Proofs}
\label{sec-proofs}
We first prove Proposition~\ref{prop:properties-I}. Afterwards, we state and prove a preliminary that we use in most of our theorems that we prove afterwards. 

\subsection*{Proof of Proposition~\ref{prop:properties-I}}
We first show convexity. Let $b_1, b_2$ be such that $I(b_j) \geq 1, j=1,2$. Let $a\in (0,1)$. We have to show that 
$I(ab_1+(1-a) b_2)\geq 1$. 
Since $\lambda$, defined in~\eqref{assumption-lighttails}, is continuous, non-decreasing, and satisfies the relation $\lambda (rx) = r^\beta \lambda (x)$, there exists an $x_{a}\geq 0$ such that 
$I(ab_1+(1-a) b_2) = \lambda(x^{(a)})$. Moreover, by definition of $I$ in \eqref{def-I}, there exists a sequence $(y_\varepsilon)_{\varepsilon>0}$ such that $(ab_1+(1-a) b_2)^Ty_\varepsilon\ge 1$ and $\lambda(y_\varepsilon)\downarrow \lambda(x^{(a)})$ as $\varepsilon\downarrow 0$. For each $\varepsilon>0$, 
$\max_{j=1,2}b_j^T y_\varepsilon\ge 1$. Suppose without loss of generality that $b_j^T x^{(a)} \geq 1$. Then, for some $\delta_\varepsilon$ tending to $0$ as $\varepsilon\downarrow 0$,
\[
   I(ab_1+(1-a) b_2) \ge \lambda (x^{(a}) -\delta_\varepsilon \geq 
   \inf_{x: b_1^T x \geq 1}  \lambda (x)-\delta_\varepsilon = I(b^{(1)}) -\delta_\varepsilon \geq 1-\delta_\varepsilon.
\]
Since this is true for each $\varepsilon>0$, $I(ab^{(1)}+(1-a) b^{(2)})\ge 1$, which implies the desired convexity property.

\smallskip
Next, we show that $I(\cdot)$ is continuous, which implies that $\{b: I(b)\ge 1\}$ is closed.
We apply Theorem 1.17(c) of \cite{RockWets98}. Following their notation, we write $f(x,u)=\lambda (x)$ if $u^Tx \geq 1$ and $f(x,u)=\infty$ otherwise. This function is lower-semicontinuous. In addition, given $\bar u\geq 0$ and $M\in [0,\infty)$, there exists a neighborhood $V$ of $\bar u$ such that $\{(x,u): u\in V, f(x,u)\leq M\}$ is bounded, since $\lambda(x) = |x|^\beta \lambda (x/|x|)$ and $\inf_{x: |x|=1} \lambda(x)>0$ by~\eqref{eq:lambda-scaling}. Thus, $f$ is uniformly level bounded conform Definition 1.16 in~\cite{RockWets98}. Let $p(u):=\inf_xf(x,u)=\inf\{\lambda(x): u^Tx\ge 1, x\in[0,\infty)^m\}=I(u)$. By continuity of $\lambda$, the infimum is attained at a (possibly non-unique) value $\bar x$ for a given $u$. 
Since the definition of $\bar x$ implies that $ u^T \bar x\geq 1$, we see that $f$ is continuous at $(u, \bar x)$, so that all conditions of Theorem 1.17(c) of \cite{RockWets98} hold. We conclude that $p(u)=I(u)$ is continuous at $u\geq 0$. Therefore $\{b: I(b)\ge 1\}$ is closed. Boundedness of $\{b: I(b)\ge 1\}$ follows from the relation $\lambda (rx) = r^\beta \lambda (x)$, $\beta>0$, and thus the set is compact.

\smallskip 
We now turn to a proof of \eqref{tail bl}. We first construct a lower bound. Given $\varepsilon>0$, let $x$ be such that $b^Tx\geq 1$ and $\lambda(x) \leq I(b)+\varepsilon$. Observe that if $L>rx$, then $b^TL>rb^Tx\ge 1$ for any $b\ge 0$, which implies $\{b^TL>r\}\supseteq \{L>rx\}$.
Using Assumption (\ref{assumption-lighttails}), it follows that, as $r\to\infty$
\begin{equation}\label{eq:btl-lower}
    {\bf P}( b^T L > r) 
    \geq  {\bf P}( L > rx) = e^{-q(r) \lambda (x)(1+o(1))} \geq e^{-q(r) (I(b)+\varepsilon)(1+o(1))}.
\end{equation}
Since this holds for every $\varepsilon>0$, the lower bound of \eqref{tail bl} follows.

\smallskip 
We turn to the upper bound. If $b$ contains elements with zeroes, we may leave out the corresponding coordinates of $L$, and the reduced random vector $\tilde L$ still satisfies Assumption 
\eqref{assumption-lighttails}, with function $\lambda ((0,x))$ (with zeroes at the corresponding coordinates of $b$). So, we may assume without loss of generality that the coordinates of $b$ are all strictly positive. 
As a result,  given $\varepsilon >0$, the set $A_{\varepsilon}=\{x\geq 0: b^Tx=1- \varepsilon \}$ is compact. 
Let  $\{x_1, \ldots,x_N\}\subset A_{\varepsilon}$, for some large constant $N=N(\varepsilon)$, be constructed such that 
\[
\{ y\ge 0: b^Ty>1\}\,\subseteq \,\bigcup_{j\in[N]}\{ y>x_j: b^Ty>1\}\,\,\,\Big(\subseteq \,\big\{y>x_j\text{ for some }j\in[N]\big\}\Big),
\]
which is possible when $N$ is sufficiently large.
  We obtain by a union bound and~\eqref{assumption-lighttails}, as $r\to\infty$,
\begin{align*}
     {\bf P}\big( b^T (L/r) > 1) \leq \sum_{j=1}^{N}  {\bf P}( L/r > x_j )\leq N \max_j {\bf P}( L > rx_j )&=  N e^{-q(r) (1+o(1)) \min_j \lambda (x_j)}.
\end{align*}
Since all $x_j$ satisfy $b^Tx_j=1-\varepsilon$, we have $\lambda (x_j) \geq I(b/(1- \varepsilon)) = (1-\varepsilon)^\beta I(b)$. Thus, as $r\to\infty$,
\begin{equation}\nonumber
    {\bf P}( b^TL> r) \leq 
    e^{-q(r) (1+o(1)) (1-\varepsilon)^\beta I(b)},
\end{equation}
providing a matching upper bound to~\eqref{eq:btl-lower} and thus proving~\eqref{tail bl}.
Lastly, we use the principle of the largest term (see e.g.\ \cite{ganesh2004}, Lemma 2.1) to conclude the asymptotics in \eqref{light-tail-phi}.
\qed

\medskip
We proceed to a preliminary result that is used in the convergence results of Theorems~\ref{thm:lt-ccp},~\ref{thm:lt-scen},~\ref{thm:ht-ccp}, and~\ref{thm:ht-cvar}. 
There are related results in the literature, cf.\  \cite{Bonnans}, \cite{RockWets98}, but our direct approach seems the most efficient for our purposes. 
For $x,y\in\mathbb{R}^m$, we say that $x\le y$ if $x_i\le y_i$ for all $i\in [m]$. If one of these inequalities is strict, we write $x<y$. We call $f:\mathbb{R}^m\mapsto\mathbb{R}$ (strictly) increasing if $x\le y$ (resp.\ $x<y$) implies that $f(x)\le f(y)$ (resp.\ $f(x)<f(y)$).

\begin{proposition}\label{prop:1}
    Let $g_\ell, g: [0,\infty)^m\mapsto \mathbb{R}$ be functions, such that $y\mapsto g_\ell( y)$ is non-decreasing for each $\ell$, and $\eta\mapsto g(\eta y)$ is strictly increasing for each $y$ as a function of the scalar $\eta\in [0,\infty)$. Let $c\in[0,\infty)^m$, and define $v_\ell:=\sup \{c^Ty: g_\ell(y)\le 1, y\ge 0\}$ and $v_\ast := \sup \{c^Ty: g(y)\le 1, y\ge 0\}$. Then, $v_\ell\to v_\ast$.
\end{proposition}
\begin{proof}
    Assume first that $v_\ast=\infty$. For each fixed $a>0$, there exists $y^{(a)}\ge0$ such that $c^Ty^{(a)} \geq a$ and $g(y^{(a)})\leq 1$.
    Fix $\eta\in(0,1)$. By the pointwise convergence $g_\ell(\eta y^{(a)})\to g(\eta y^{(a)})<1$, where the strict inequality follows from the fact that $g(\eta y^{(a)})$ is strictly increasing in $\eta$. In other words, $\eta y^{(a)}$ is feasible for $g_\ell$ for all sufficiently large $\ell$. Moreover, for such $\ell, v_\ell\ge c^T \eta y^{(a)}=\eta c^Ty^{(a)}\ge \eta a$. Since $a$ is arbitrary, we obtain that $\liminf_{\ell\to\infty} v_\ell =\infty=v_\ast$.
    Assume now that $v_\ast<\infty$. Analogously to above, $v_\ell\ge \eta v_\ast$ for each $\eta\in(0,1)$ and all $\ell$ sufficiently large depending on $\eta$. As a result, $\liminf_{\ell \to\infty} v_\ell \ge v_\ast$. 
    We  turn to the upper bound. 

    Let us partition $c=(\bar{c},0)$, where the part $\bar{c}$ has strictly positive elements. We then note that $v_* = \sup\{\bar{c}^T \bar{y} : g(\bar{y},0) \leq 1, \bar{y} \geq 0\}$. Because $g_\ell$ is  non-decreasing 
    we have that 
$v_\ell = \sup\{\bar{c}^T \bar{y}: g_\ell(\bar{y},0) \leq 1, \bar{y}\geq 0\}.$ 
 So, we may assume without loss of generality that $c$ has strictly positive elements.
 
    Assume for contradiction that there exists $\varepsilon>0$ such that $(v_\ell)_{\ell\ge 1}$ contains a subsequence such that $v_{\ell_i} \ge v_\ast+3\varepsilon$ for all $i$. Since $v_{\ell_i}=\sup \{c^Ty: g_{\ell_i}(y)\le 1, y\ge 0\}$,  there must exist $y_{\ell_i}\in [0,\infty)^m$ such that $g_{\ell_i}(y_{\ell_i})\le 1$ and $c^T y_{\ell_i}\ge v_\ast+2\varepsilon$. For each $i$, there exists $\eta_i\leq 1$ and $y'_{\ell_i}\in[0,\infty)^m$ such that $y_{\ell_i}'=\eta_i y_{\ell_i},  g(y'_{\ell_i}) = g(\eta_i y_{\ell_i}) \le g(y_{\ell_i}) \le 1$ and $c^Ty'_{\ell_i}=v_\ast+2\varepsilon$.

    Let $A:=\{y\in[0,\infty)^n: c^Ty=v_\ast + 2\varepsilon\}$. As all elements of $c$ are strictly positive, $A$ is a compact set. 
    Define $f(y):=(\lfloor y_1\cdot |c|/\varepsilon\rfloor, \ldots, \lfloor y_n\cdot |c|/\varepsilon\rfloor)\cdot \varepsilon/|c|$, i.e., each coordinate is rounded down to the nearest multiple of $\varepsilon/|c|$. Define $\hat A:=\{f(y): y\in A\}\ni f(y'_{\ell_i})=: \hat y_{\ell_i}$. Since $A$ is compact, $\hat A$ is a finite set, and thus we can extract a further subsequence such that $\hat y_{\ell_{i_j}}$ is identical to some $\hat y\in \hat A$ along the subsequence. By the rounding procedure, $c^T \hat y\ge v_\ast+\varepsilon$, and since $g_{\ell_{i_j}}$ is non-decreasing, $g_{\ell_{i_j}}(\hat y)\le 1$. 
        By the pointwise convergence, $g(\hat y)=\lim_{j\to\infty}g_{\ell_{i_j}}(\hat y)\le 1$, and thus $\hat y$ is a feasible solution for $g$. This contradicts that $v_\ast$ is the supremum of $\{c^Ty: g(y)\le 1, y\ge 0\}$. As a result, $\limsup_{\ell\to\infty}v_\ell\le v_\ast+3\varepsilon$ for any $\varepsilon>0$. Thus, $\limsup_{\ell\to\infty}v_\ell\le v_\ast$. 
\end{proof}

\subsection{Proof of Theorem~\ref{thm:lt-ccp}}\label{sec:proof-1}

\smallskip
We now turn to the proof of Theorem~\ref{thm:lt-ccp}.
We aim to apply Proposition 1 with $\ell=1/\delta$, 
$g_\ell(y) = (\log \delta)/ \log {\bf P}(\phi(y,L) > q^{-1}(\log 1/\delta))$, and $g(y)=1/J(y)$. The pointwise convergence $g_\ell(y)\rightarrow g(y)$ follows from~\eqref{tail bl}.
The fact that $\eta\mapsto g(\eta y)$ is strictly increasing for all $y$ follows from the homogeneity property of $J$ implied by \eqref{def-I}.
Moreover, $\eta\mapsto g_\ell(\eta y)$ is non-decreasing because $\eta \mapsto \phi(\eta y, \ell)$ is non-decreasing for each $y$ and $\ell$. Thus, by Proposition 1, when $\delta$ is sufficiently small,
\[
\begin{aligned}
q^{-1}(\log1/\delta)\cdot v_\delta&=q^{-1}(\log1/\delta)\cdot\sup\Big\{c^Tx: {\bf P}\big(\phi(x, L)>1\big)\le \delta, x\in X\Big\} \\&=
\sup\Big\{c^Ty: \delta\ge {\bf P}\big(\phi(y,L)>q^{-1}(\log 1/\delta)\big), y\ge 0\Big\}\\
&=\sup\big\{c^Ty: g_\ell(y)\le 1, y\ge 0\big\} \\
&\to \sup\big\{c^Ty: g(y)\le 1, y\ge 0\big\}=\sup\big\{c^Ty: J(y)\ge 1, y\ge 0\big\}=c^Ty^\ast.
\end{aligned}
\]
In the second equality above, we used that the constraint $x\in X$ is no longer relevant in (\ref{limitingoptimizationproblem-lt}) as $r\rightarrow\infty$, since we assumed $[0,h]^m \subseteq X$ for some $h>0$, and $x=y/r$.
  \qed

\subsection{Proof of Theorem~\ref{thm:lt-cvar}}\label{sec-proof2}

Since $v_{{\bf CVaR}, \delta}\leq v_\delta$, it suffices to prove an asymptotic lower bound. For this, we utilize Theorem~\ref{thm:lt-ccp}. Let $y^*$ be a solution of \eqref{limitingoptimizationproblem-lt} and for fixed $\gamma>0$, define $y := (1-2\gamma) y^*$. Moreover, let  $r=r_{\rm{LT}, \delta}:=q^{-1}(\log 1/\delta)$. We will prove that $x=y/r_{\rm{LT}, \delta}$ will be feasible for \eqref{cvar-opt}, for $\delta$ sufficiently small. 
We can write (\ref{cvar}) with $x=y/r_{\rm{LT}, \delta}$ as
\begin{equation}\label{eq:pr2-cvar}
 {\bf CVaR} (\phi(y/r_{\rm{LT}, \delta},L)) =  \inf_{\tau\in\mathbb{R}} \left[\tau + \frac 1\delta \int_\tau^\infty {\bf P} (\phi (y,L) > ru) du \right]\leq 1-\gamma + \frac 1\delta \int_{1-\gamma}^\infty {\bf P} (\phi (y,L) > ru) du.
\end{equation}
Using \eqref{tail bl}, for every $\varepsilon >0$ there exists an $r^{(1)}_\varepsilon$ such that
\begin{equation}
     {\bf P} (\phi (y,L) > ru) = e^{-q(ru) J(y)(1+o(1))} \leq  e^{-q(ru) J(y) (1-\varepsilon)}, \qquad ru \geq r^*_\varepsilon.
\end{equation}
In addition, 
due to the Potter bounds (see e.g.\ \cite{Bingham_Goldie_Teugels_1987}, Theorem~\ref{thm:lt-ccp}.5.6), there exists a constant $r^{(2)}_{\varepsilon}$ such that, for $u\geq 1/2$, and $r\geq r^{(2)}_{\varepsilon}$,
\begin{equation}
     q(ru)/ q(r)  \geq (1-\varepsilon) u^{\beta-\varepsilon}.
\end{equation}
Consequently, for $r$ sufficiently large and $u>1-\gamma$, using also 
that $J(y) = (1-2\gamma)^{-\beta} J(y^*)= (1-2\gamma)^{-\beta}$,
\begin{equation}\label{eq:pr2-prob}
     {\bf P} (\phi (y,L) > ru)  \leq  e^{-q(r) u^{\beta - \varepsilon} (1-\varepsilon)^2/(1-2\gamma)^\beta}.
\end{equation}
Observe $\delta = \exp(-q(r))$, so that we get from \eqref{eq:pr2-prob} and \eqref{eq:pr2-cvar} 
\begin{equation}
 {\bf CVaR} (\phi(y/r_{\rm{LT}, \delta},L)) \leq 1-\gamma + \int_{1-\gamma}^\infty e^{ q(r)}  e^{-q(r) u^{\beta -\varepsilon}(1- \varepsilon)^2/(1-2\gamma)^\beta} du
 \end{equation}
Now, given $\gamma$, choose $\varepsilon$ small enough such that $(1-\gamma)^{\beta - \varepsilon} (1-\epsilon)^2/(1-2\gamma)^\beta>1$.
In that case, the RHS of the last display converges to $1-\gamma$ as $r\to\infty$. We conclude that $x=(1-2\gamma) y^*/r$ is feasible for $\delta$ sufficiently small. Since this is true for each $\gamma>0$, the proof is complete. 
\qed

\subsection{Proof of Theorem~\ref{thm:lt-scen}}
Let $L^{(j)},j=1,...,k$ be $k$ i.i.d.\ samples of $L$.
Combining~\eqref{tail bl} with Theorem~\ref{thm:lt-ccp} of
\cite{Resnick1973}, 
 we have that
\begin{equation}
\label{evtlight}
    \frac 1{q^{-1} (\log k)} \max_{j=1,...,k} b^T L^{(j)}\stackrel{a.s.}{\rightarrow} \frac  1{I(b)^{1/\beta}}, \qquad k\to\infty.
\end{equation}
almost surely.
We apply \eqref{evtlight} to establish the following pointwise convergence for $y\ge 0$:
\begin{align}
    g_k(y) &:=\frac 1{q^{-1} (\log k)} \max_{j=1,...,k}  \max_{i=1,...,d}y^TA^{(i)} L^{(j)}
    \stackrel{a.s.}\rightarrow \max_{i=1,...,d} 1/I(y^T A^{(i)})^{1/\beta}
    =: g(y).
\end{align}
The functions $g_k(y)$ are non-decreasing for all $k$. We claim that $\eta \mapsto g(\eta y)$ is strictly increasing, which is implied if $I(y^T A^{(i)})$ is strictly decreasing for all $i=1,\ldots, d$. The latter is implied by the equality in the definition of $I$ in~\eqref{def-I}.
The statement for arbitrary $k$ follows now immediately from Proposition~\ref{prop:1}. 

Next, assume $k=k(\delta, \varepsilon)$  defined in \eqref{cc-choice}. Then, as $\delta\downarrow 0$, for any $\varepsilon=\varepsilon(\delta)$ satisfying $\delta\log(1/\varepsilon)\to0$,
\begin{align}
    \log k(\delta, \varepsilon)&= \log\big((2/\delta)\log(1/\varepsilon)+2m +(2m/\delta)\log(2/\delta)\big)\nonumber \\
    &= \log(1/\delta) + \log\big(2\log(1/\varepsilon)+ 2\delta/m + 2m\log(2/\delta)\big) \sim \log (1/\delta).\label{eq:asymptotics-kde}
\end{align}
 In view of Theorem \ref{thm:lt-ccp} and the first part of Theorem~\ref{thm:lt-scen}, it suffices to show that, as $\delta\downarrow 0$,
\begin{equation}\label{eq:asymptotics-kde-desired}
q^{-1}(\log 1/\delta)\sim q^{-1}(\log k(\delta, \varepsilon)),
\end{equation}
where $q(r)=h(r)r^\beta$, $\beta>1$, defined below \eqref{assumption-lighttails}.
By Theorem 1.5.12 in \cite{Bingham_Goldie_Teugels_1987}, $q^{-1}$ is regularly varying with index $1/\beta$. The desired asymptotics in \eqref{eq:asymptotics-kde-desired} follow now from \eqref{eq:asymptotics-kde}.
\qed

\subsection{Proof of Theorem~\ref{thm:ht-ccp}}\label{sec-proof4}
We recall that $r_\delta=\bar F^{-1}(\delta)$.
We first derive the pointwise limit of 
${\bf P}(\phi(y,L) > r_\delta)/{\bf P}(|L| > r_\delta)$
as $\delta\downarrow 0$,  using \eqref{rvmult}.
Let $\varepsilon \in (0, 1/\max_i y^T A^{(i)} \mathrm{e})$, with $\mathrm{e}$ a vector of ones. 
As we assumed that the entries of all $A^{(i)}$ are non-negative,  $|L| \leq \varepsilon r_\delta$ implies that $\phi(y,L) \leq r_\delta$. Consequently, 
\begin{align}
    \frac{{\bf P}(\phi(y,L) > r_\delta)}{{\bf P}(|L| > r_\delta)} = \frac{{\bf P}(\phi(y,L) > r_\delta ; |L| > \varepsilon r_\delta)}{{\bf P}(|L| > r_\delta)}\label{eq:pareto-1}&= 
     {\bf P}\big(\phi(y,L) > r_\delta \,\big|\, |L| > \varepsilon r_\delta\big)     
\frac{{\bf P}(|L| > \varepsilon r_\delta)} {{\bf P}(|L| > r_\delta)}\\
  &= 
     {\bf P}\Big(\phi\Big(y,\frac{L}{|L|}\Big) \frac{|L|}{r_\delta} >1 \,\big|\, |L| > \varepsilon r_\delta\Big)     
\frac{{\bf P}(|L| > \varepsilon r_\delta)} {{\bf P}(|L| > r_\delta)}.\label{eq:pareto-3}
\end{align}
Due to univariate regular variation, the ratio converges to 
$\varepsilon^{-\alpha}$. Due to the multivariate regular variation property \eqref{rvmult}, the conditional (on $ |L| > \varepsilon r_\delta$) distributional limit of $(L/|L|, |L|/r_\delta)$ 
is $(\Theta, \varepsilon T)$, with $T$ a random variable with density $\alpha u^{-1-\alpha}$ independent of $\Theta$, see below \eqref{rvmult} and \cite{basrak2002characterization}.
Consequently, 
\begin{align}
 \frac{{\bf P}(\phi(y,L) > r_\delta)}{{\bf P}(|L| > r_\delta)} \rightarrow \varepsilon ^{-\alpha}   
 {\bf P}(\phi(y,\Theta) \varepsilon T >1) 
 &= \varepsilon ^{-\alpha} \int_1^\infty \alpha u^{-1-\alpha} 
  {\bf P}(\phi(y,\Theta) \varepsilon  >1/u)\rd u\\
  &= \int_0^{1/\varepsilon} \alpha v^{\alpha-1} 
  {\bf P}(\phi(y,\Theta)  >v)\rd v\\
  &= \int_0^{\infty} \alpha v^{\alpha-1} 
  {\bf P}(\phi(y,\Theta)  >v)\rd v= {\bf E} [ \phi(y,\Theta)^\alpha].\label{eq:pareto-2}
\end{align}
In the derivation above, the second equality follows from
the change of variable $v=1/\varepsilon u$, the third equality comes from the fact that 
$\phi(y,\Theta) \leq \max_i y^T A^{(i)} e \leq 1/\varepsilon$, and the last equality from a representation for the $\alpha$th moment of a non-negative random variable. 

To finish the proof, we define $k=1/\delta$, $g_k(y)={\bf P}(\phi(y,L) > r_\delta)/{\bf P}(|L| > r_\delta)$ and $g(y)={\bf E}[\phi(y,\Theta)^\alpha]$. By monotonicity of $y\mapsto\phi(y, L)$, these functions are all increasing.  By Proposition~\ref{prop:1}, when $\delta$ is  small,
\begin{align*}
    r_\delta v_\delta = r_\delta \sup\Big\{c^Tx: {\bf P}\big(\phi(x, L)>1\big)\le \delta, x\in X\Big\}
    &=\sup\Big\{c^Ty: {\bf P}\big(\phi(y, L)>r_\delta\big)/{\bf P}(|L|>r_\delta)\le 1, y\ge 0\Big\} \\
    &\to \sup\Big\{c^Ty: {\bf E}\big[\phi(y,\Theta)^\alpha\big]\le 1, y\ge 0\Big\}=c^Ty^\ast.
\end{align*}
We again used that the constraint $x\in X$ is replaced by $y\ge 0$, as we assumed $[0, h]^m \subseteq X$ for some $h>0$, and $x=y/r$.\qed

\subsection{Proof of Theorem~\ref{thm:ht-cvar}}
Similar to~\eqref{eq:pr2-cvar}, we write (\ref{cvar}) with $x=y/r_\delta$ as
\begin{equation}
 {\bf CVaR} (\phi(y/r_\delta,L)) =  \inf_{\tau\in\mathbb{R}}\left\{ \tau + \frac 1\delta \int_\tau^\infty {\bf P} (\phi (y,L) > r_\delta u) du\right\} = \tau_\delta + 
 \frac 1{r_\delta\delta} \int_{r_\delta\tau_\delta}^\infty {\bf P} (\phi (y,L) > v) dv,
\end{equation}
where $\tau_\delta= \inf \{t:  {\bf P}(\phi(y,L) > r_\delta t) \leq \delta\}$ attains the infimum.
By the reasoning in the proof of Theorem \ref{thm:ht-ccp} from \eqref{eq:pareto-1} until \eqref{eq:pareto-2} with $y$ replaced by $y/t$, we obtain
\begin{equation}
    {\bf P}(\phi(y,L) > r_\delta t)\sim {\bf E} [ \phi(y/t,L)^\alpha] \bar F (r_\delta) \sim {\bf E} [ \phi(y,L)^\alpha] t^{-\alpha} \delta.
\end{equation}
Consequently, 
    $\tau_\delta \sim {\bf E} [ \phi(y,L)^\alpha]^{1/\alpha}$.
Since $r_\delta\tau_\delta\rightarrow\infty$ and ${\bf P} (\phi (y,L) > v)$ is regularly varying with index $-\alpha<-1$ by the homogeneity of $\phi$, we can apply Karamata's theorem to obtain 
\begin{equation}
    {\bf CVaR} (\phi(y/r_\delta,L)) \sim \tau_\delta + \frac{r_\delta\tau_\delta}{r_\delta\delta} \frac{1}{\alpha-1} {\bf P}(\phi(y,L) > r_\delta\tau_\delta) \sim \frac{\alpha}{\alpha-1} \tau_\delta \sim \frac{\alpha}{\alpha-1} {\bf E} [ \phi(y,L)^\alpha]^{1/\alpha}.
\end{equation}
Now, we apply Proposition \ref{prop:1} with $k=1/\delta$, $g_k(y) = {\bf CVaR} (\phi(y/r_\delta,L))$, and $g(y) = \frac{\alpha}{\alpha-1} {\bf E} [ \phi(y,L)^\alpha]^{1/\alpha}$. From the variational representation of $ {\bf CVaR} (\phi(y/r_\delta,L))$ and the fact that $\phi$ is nondecreasing in $y$ for every fixed $L$, it follows that $g_k$ is nondecreasing in $y$, and $g(ay)$ is strictly increasing in $a$ due to the scaling property of $\phi$. Since we also have established that $g_k(y)\rightarrow g(y)$, we can apply Proposition 1 to conclude that the normalized version of \eqref{cvar-opt} converges to the optimal value of 
\begin{equation}
\label{limitingcvarproblem}
\mbox{maximize } c^Ty \hspace{1cm} {\mbox{ subject to }}   y\geq 0 \mbox{ and } {\bf E}[\phi(y,\Theta)^\alpha]\leq (1-1/\alpha)^\alpha.
\end{equation}
If $y^\ast$ solves \eqref{limitingoptimizationproblem}, then $(1-1/\alpha)y^\ast$ is an optimal solution of \eqref{limitingcvarproblem} by the homogeneity of $\phi$. Combining this with Theorem~\ref{thm:ht-ccp} completes the proof. 
\qed

\subsection{Proof of Theorem~\ref{thm:ht-scen}}\label{sec-proof6}
Write $a_k:= \bar F^{-1} (1/k)$ for the normalizing constants in~\eqref{sample-problem-scaled}. To show convergence in distribution of the optimal value $V_{SCP}(k)$ to $V^\ast$, it suffices to show that the non-feasible set spanned by the constraints $\phi(y, L^{(i)}/a_k)\le 1$, $i=1,\ldots, k$, converges in distribution to the non-feasible set spanned by the constraints $\phi(y, z)\le 1$, $z\in N$, where $N$ is the Poisson point process defined above~\eqref{eq:pointprocess}. To do so, 
we rely on the theory of random sets, as developed in \cite{Molchanov}. 

To translate our setting to that of Section 4.2 in \cite{Molchanov}, define the multi-function $M(L) := \{ y: \phi (y,L) \ge 1\}$.
By the homogeneity of $\phi$, the closure of the set of non-feasible solutions for \eqref{sample-problem-scaled}
can be written as 
\begin{equation}
    U_k = \bigcup_{i=1}^k \Big\{ y: \phi\Big(y,\frac{L^{(i)}}{a_k}\Big) \ge 1\Big\} = \frac{1}{a_k}  \bigcup_{i=1}^k M(L^{(i)}),
\end{equation} 
where we use the notation $bS = \{x: b^Tx \in S\}$. 
Applying Theorems~\ref{thm:ht-ccp}.2.9 and 4.2.10 of \cite{Molchanov}, we obtain that $U_k$ converges in distribution to $U = \{y: \phi (y,z) \ge 1 \mbox{ for some } z\in N\}$. 
We define for $t>0$ the half-space  $K_t =\{ y: c^Ty > t\}$. 
 Therefore, we see that for $t>0$,
\begin{align}
    {\bf P} (V_{SCP}(k) \leq t) &={\bf P}\big(\forall y\in K_t \exists i\in\{1,\ldots, k\}:  \phi\big(y, L^{(i)}/a_k\big)> 1\big) \\
    &={\bf P}\big(\forall y\in K_t \exists i\in\{1,\ldots, k\}:  \phi\big(y, L^{(i)}/a_k\big)\ge 1\big)\\
    &= {\bf P}(U_k \cap K_t \neq \emptyset) \rightarrow 
    {\bf P}(U \cap K_t \neq \emptyset) = {\bf P} (V^* \leq t). 
\end{align}
This completes the proof of Theorem~\ref{thm:ht-scen}. \qed

\end{document}